\documentclass[11pt]{article}

\usepackage{amsmath}
\usepackage{amssymb}
\usepackage{eucal}

\begin{document}

\baselineskip 16pt

\title{On  Hall   subgroups of a finite  group }

\author{Wenbin Guo\thanks{Research of the first author is supported by
a NNSF grant of China (Grant \# 11071229)  and Wu Wen-Tsun Key Laboratory of Mathematics, USTC, Chinese Academy of Sciences.}\\
{\small Department of Mathematics, University of Science and
Technology of China,}\\ {\small Hefei 230026, P. R. China}\\
{\small E-mail:
wbguo@ustc.edu.cn}\\ \\
{ Alexander  N. Skiba \thanks{Research of the second author
 supported by Chinese Academy of Sciences Visiting Professorship for Senior
 International Scientists (grant No. 2010T2J12)}
}\\
{\small Department of Mathematics,  Francisk Skorina Gomel State University,}\\
{\small Gomel 246019, Belarus}\\
{\small E-mail: alexander.skiba49@gmail.com}}

\maketitle

\begin{abstract}

In the paper  new  criteria of existence and conjugacy of Hall
subgroups of finite groups are given.

\end{abstract}

\let\thefootnoteorig\thefootnote
\renewcommand{\thefootnote}{\empty}


\let\thefootnote\thefootnoteorig

\renewcommand{\thefootnote}{\empty}

\section{Introduction}

Throughout this paper, all groups are finite and $G$ always denotes
a finite group. We use  $\Bbb{P}$ to denote the set
of all primes, $\pi$ is a  non-empty subset of
$\Bbb{P}$ and $\pi'=\Bbb{P}  \setminus \pi$. 
 A natural number $n$   is called a  \emph{$\pi$-number} if  $(n, p)=1$ for any prime 
$p\in \pi'$.  If $|G|$ is  a  $\pi$-number, then $G$ is said to be  a  
\emph{$\pi$-group}.

 Let  $A$ and $B$ be  subgroups  of $G$. Then  $A$
 is said to permute with  $B$ if  $AB=BA$.    If $A$ permutes with all 
(Sylow) subgroups of $B$, then $A$ is called \emph{quasinormal}  
(\emph{$S$-quasinormal}, respectively) in $B$.

A group $G$ is said to be:

(a)  an \emph{$E_{\pi}$-group} 
   if $G$ has a  Hall $\pi$-subgroup;

(b)  a \emph{$C_{\pi}$-group}    if $G$ is an $E_{\pi}$-group   and any
two Hall $\pi$-subgroups of $G$ are conjugate;

(c)  a \emph{$D_{\pi}$-group}   
  if $G$ is an $E_{\pi}$-group  and any
$\pi$-subgroup of $G$ is contained in some Hall $\pi$-subgroup of $G$.

A group $G$ is said to be \emph{ $\pi$-separable} if $G$ has a  chief 
 series  $$1=G _0\leq H_1 \leq \ldots \leq H_{t-1}\leq H_t =G,
\eqno(*)$$ where each index $|H_{i}:H_{i-1}|$ is either a
$\pi$-number or a $\pi'$-number.

 The most important result of the  theory of  $\pi$-separable  groups 
 is the following classical result.

{\bf Theorem A}  (P. Hall, S.A. $\breve{C}$unihin).  
  {\sl Any  $\pi$-separable group is  a     $D_{\pi}$-group.}

Our main goal here is to prove the following generalization of this 
theorem.

 {\bf Theorem 1.1.} {\sl Suppose that $G$ has a   subgroup series  $$1=H _0 < H_1 \leq
\ldots \leq H_{t-1}\leq H_t =G,$$  where $|H_{i+1}:H_i|$ is either a $\pi$-number
 or a $\pi'$-number for all $i=1, \ldots , t$.  Suppose also that $G$ has a 
 subgroup $T$  such that $H_{1}T=G$  and  $|G:T|$ is  a $\pi'$-number.   
   If  $H_{i}$ is quasinormal in $T$,  for    all $i= 1,  \ldots , t$,
 then $G$ is  a $D_{\pi}$-group.}

The following theorem  is one of the main steps in the proof of Theorem 
1.1.

 {\bf Theorem 1.2.} {\sl Suppose that $G$ has a   subgroup series  $$1=H _0 < H_1 \leq
\ldots \leq H_{t-1}\leq H_t =G,$$  where $|H_{i+1}:H_i|$ is either a $\pi$-number
 or a $\pi'$-number for all $i=1, \ldots , t$.  Suppose also that $G$ has a 
 subgroup $T$  such that $H_{1}T=G$  and  $|G:T|$ is  a $\pi'$-number.   
  If $H_{i}$ is $S$-quasinormal in  
$T$,  for    all $i= 1,  \ldots , t$, then $G$ is 
 a $C_{\pi}$-group.}

Theorem 1.1  strengthens  the main result in \cite{Proc}.
  Note also that  Example 1.1 in \cite{Proc} shows that,
 under the conditions in Theorems 1.1 or 1.2
 the group $G$ is not necessary $\pi$-separable.

\section{Premilaries}

{\bf Lemma 2.1.}
{\sl Let $N$ be a normal $C_{\pi}$-subgroup of $G$.}

(i) {\sl  If  $G/N$ is a
$C_{\pi}$-group, then $G$ is a $C_{\pi}$-group}   (S. A. $\breve{C}$unihin 
\cite{ChuII}).

(ii) {\sl  If  $G/N$ is an
$E_{\pi}$-group, then $G$ is an $E_{\pi}$-group} 

(iii) {\sl  If  $G$ has a nilpotent    Hall 
$\pi$-subgroup, then $G$ is a $D_{\pi}$-group}   (\cite{Wi}).

(iv) {\sl  If  $G$ has a     Hall 
$\pi$-subgroup with cyclic Sylow subgroups, then $G$ is a $D_{\pi}$-group} 
(S. A. Rusakov  \cite{Rus}).

{\bf Lemma 2.2.}
{\sl Let $N$ be a normal $C_{\pi}$-subgroup of $G$  and $N_{\pi}$ a Hall 
$\pi$-subgroup of $N$.}

(i) {\sl  If  $G$ is a
$C_{\pi}$-group, then $G/N$ is a $C_{\pi}$-group}  (See \cite[Lemma 
9]{HallIII}).

(ii) {\sl If every Sylow subgroup of $N_{\pi}$  is cyclic and  $G/N$ is a 
  $D_{\pi}$-group, then $G$ is a  $D_{\pi}$-group }   
 (See \cite{DAN} or \cite[Chapter IV, Theorem 18.17]{26}).

(iii) {\sl    $G$ is a
$D_{\pi}$-group if and only if  $G/N$ is a  $D_{\pi}$-group}  (See 
\cite{Rev}).

{\bf Lemma 2.3} (O. Kegel \cite{KegI}). {\sl Let $A$ and
$B$ be subgroups of $G$ such that  $G\ne AB$ and $AB^{x}=B^{x}A$,
for all $x\in G$. Then $G$ has a proper normal subgroup $N$ such
that either $A\leq N$ or $B\leq N$.}

Let $A$ be a subgroup of $G$. A subgroup $T$ is  said to be a minimal
supplement of $A$ in $G$ if $AT=G$ but $AT_{0}\ne G$ for all proper
subgroups $T_{0}$ of $G$.

The following lemma is obvious.

{\bf Lemma 2.4.} {\sl If $N$ is normal in $G$  and $T$ is a minimal
supplement of $N$ in $G$, then $N\cap T\leq \Phi (T)$.}

{\bf Lemma 2.5} (P. Hall \cite{Hall}). {\sl Suppose that $G$ has a
Hall $p'$-subgroup for each prime $p$ dividing $|G|$. Then $G$ is
soluble. }

Let $A$ and $B$ be subgroups of $G$ and $\varnothing \ne X\subseteq
G$. Following \cite{GSSI}, we say that $A$ is $X$-permutable (or $A$
$X$-permutes) with $B$ if $AB^{x}=B^{x}A$ for some $x\in X$.

The following lemma is also evident.

{\bf Lemma 2.6.} {\sl Let $A, B, X$ be subgroups of $G$  and $K
\trianglelefteq G$. If $A$ is $X$-permutable with $B$, then $AK/K$
is $XK/K$-permutable with $BK/K$ in $G/K$.}

{\bf Lemma 2.7}  (Kegel \cite{KegII}). {\sl If a
subgroup  $A$ of $G$ permutes with all Sylow subgroups of $G$, then
$A$  is subnormal in $G$.}

{\bf Lemma 2.8} (H. Wielandt \cite{Wiel}). {\sl  If a $\pi$-subgroup
$A$ of $G$ is subnormal in $G$, then $A\leq O_{\pi}(G)$.}

{\bf Lemma 2.9} (V. N. Knyagina and V. S. Monakhov \cite{mon}). {\sl
  Let $H$, $K$  and $N$
 be 
subgroups of $G$.  If $H$ is a Hall subgroup of $G$, then $$N\cap HK=(N\cap H)(N\cap K).$$}

{\bf Lemma  2.10. } {\sl Let   $U\leq B\leq G$ and  $G=AB$, where the subgroup $A$ permutes
 with $U^{b}$ for all  $b\in B$. Then  $A$ permutes
 with $U^{x}$ for all  $x\in G$.   }

{\bf Proof.}   Since $G=AB$,  $x=ab$ for  some $a\in A$  and $b\in B$.
        Hence $$AU^{x}= AU^{ab} =Aa(U^{b})a^{-1}=a(U^{b})a^{-1}A=U^{x}A.$$

 {\bf Lemma 2.11 } \cite[Chapter I, Lemma 1.1.19]{prod}. {\sl Let 
 $A, B\leq  G$ and  $G=AB$.  Then  $G_{p}=A_{p}B_{p}$   for
 some $G_{p}\in \mathrm{Syl}_{p}(G)$,
 $A_{p}\in \mathrm{Syl}_{p}(A)$ and $B_{p}\in \mathrm{Syl}_{p}(B)$.    }

{\bf Lemma 2.12.} {\sl Let  $X$ be a normal $C_{\pi}$-subgroup of $G$.
 Suppose that $G$ has
 a subgroup series  $$1=H
_0\leq H_1 \leq \ldots \leq H_{t-1}\leq H_t =G$$ and a supplement
$T $ of $H_{1} $  
  in $G$ such that $H_{i}$  $X$-permutes with every Sylow
subgroup of $T$      for all $i= 1,  \ldots , t$. 
 If each index
$|H_{i+1}:H_{i}|$ is either a $\pi$-number or a $\pi'$-number, then 
$G$ is an $E_{\pi}$-group.  }

{\bf Proof.} 
Consider the subgroup series $$1=H _0 X/X  \leq H_1 X/X \leq \ldots \leq
H_{t-1} X/X \leq H_t X/X=G /X \eqno(1)$$  in $G/X$. If $T\leq X$, then 
$G=H_{1}T=XH_{1}$. Hence $G$ is an $E_{\pi}$-group by Lemma 2.1(2). Now 
let $T\nleq X$ and  $P/X$ be a Sylow 
$p$-subgroup of $
TX/X$. Then by Lemma 2.11 there are Sylow  $p$-subgroups $T_{p}$ of $T$ 
and $X_{p}$  of $X$ such that  $T_{p}X_{p}$ is a Sylow $p$-subgroup of 
$TX$ and $T_{p}X_{p}X=T_{p}X=P$. Hence by  Lemma 2.6, $H_{i} X/X$ 
 permutes  with every Sylow
subgroup of $T X/X$ for all $i= 1,  \ldots , t$. On the other
hand, since  
$$|H_{i+1}X/X:H_{i}X/X|=|H_{i+1}X:H_{i}X|=|H_{i+1}:H_{i}|:|X\cap
H_{i+1}:X\cap H_{i}|,$$
 every index of Series $ (1)$ is either
$\pi$-number or $\pi'$-number. Hence the assertion follows from Theorem B in \cite{Proc}.

{\bf Lemma 2.13} (see Lemma in \cite[Chapter A]{DH}). 
 {\sl
  Let $H$, $K$  and $N$
 be 
subgroups of $G$.  If $HK=KH$ and $HN=NH$, then $H\langle K, N
 \rangle =\langle K, N \rangle H$.}

\section{Base proposition}

{\bf Proposition  3.1. } {\sl Let $X$ be a normal  $C_{\pi}$-subgroup of $G$ and 
 $A$ a subgroup of $G$ such     
  that $|G:A|$ is a $\pi $-number. Suppose that 
  $A$ has a  Hall  $\pi$-subgroup $A_{0}$ such that either $A_{0}$ is nilpotent or
 every Sylow subgroup of $A_{0}$ is cyclic.  Suppose that $A$  
$X$-permutes  with every Sylow $p$-subgroup of $G$ for all primes
$p\in  \pi$ or for all primes $p\in  \pi \setminus \{ q \}$ for some prime 
  $q$  dividing $|G:A|$. Then $G$ is a $C_{\pi}$-group.   }

{\bf Proof.}  Assume  that this proposition   is false and let $G$ be a
counterexample of minimal order. Then $|\pi \cap \pi (G)| > 1$. 

(1) {\sl  $G/R$ is a  $C_{\pi}$-group for any non-identity  normal subgroup $R$ of $G$.  }

In  order to prove this assertion, in view of the choice of $G$,  it is enough to
 show that the hypothesis is   still true for 
$(G/R, AR/R, XR/R$).   First note that $|G/R:AR/R|=|G:AR|$ is a $\pi $-number, and  $A_{0}R/R$ 
is a Hall $\pi$-subgroup of $AR/R$  since   
 $$|AR/R: A_{0}R/R|=|AR: A_{0}R|=|A:A\cap A_{0}R|=|A:A_{0}(A\cap 
R)|.$$  On the other hand, $XR/R\simeq X/X\cap R$ is a    $C_{\pi}$-group 
by Lemma 2.2 (i),  and  either $A_{0}R/R\simeq A_{0}/R\cap 
A_{0}$ is    nilpotent     or every Sylow subgroup of  $A_{0}R/R$ is 
cyclic. Finally, let 
 $P/R$ be  a  
Sylow $p$-subgroup of $G/R$, where $p\in  \pi \setminus q$. Then for some Sylow $p$-subgroup
 $G_{p}$ we have $G_{p}R/R=P/R$. Hence  $AR/R$ 
$XR/R$-permutes with  $P/R$ by Lemma 2.6.  Therefore the hypothesis holds 
for $(G/R, AR/R, XR/R)$.

(2)  {\sl $X=1$}.

Indeed, if  $X\ne 1 $, then $G/X$ is a  $C_{\pi}$-group by (1). Hence $G$ is  
 $C_{\pi}$-group  by Lemma 2.1 (i), a contradiction. 

(3) {\sl $G$ has a proper non-identity normal subgroup $N$}.

  Let $p\in 
\pi \cap \pi (G)$, where $p\ne q$. Let $P$ be a Sylow $p$-subgroup of $G$. 
 First assume that $AP=G$. Since $|\pi \cap \pi (G)| > 1$, there is a 
prime $r\in \pi \cap \pi (G)$ such that $r\ne p$, so $r$ does not divide   
       $|G:A|$. Let $R$ be a Sylow $r$-subgroup of $G$. Then for any $x\in 
G$ we have $AR^{x}=R^{x}A$. Hence $R\leq A_{G}$.  Since $G$ is not a  
$C_{\pi}$-group, $A\ne G$ by Lemma 2.1 (iii)(iv). Hence   $1\ne  A_{G}\ne G$. 
Now suppose that $AP\ne G$.  
By (2),  $P^{x}A=AP^{x}$ for all $x\in G$. Hence  we have (3)  by Lemma 2.3.

(4) {\sl $N$ is a  $C_{\pi}$-group.}

In view of the choice of $G$ it is enough to prove that the hypothesis holds
 for $(N, A_{1})$,
where $A_{1}=A\cap N$.  Since $|AN:A|=|N:A\cap N|$, $|N:A_{1}|$  is a  
$\pi $-number.  On the other hand,  $A_{0}\cap N$  is a  Hall ${\pi}$-subgroup 
  of $N$ since $$|A\cap N:A_{0}\cap N|=|A_{0}(A\cap N):A_{0}|.$$ 
It clear also that either  $A_{0}\cap N$ is nilpotent or every Sylow 
subgroup of $A_{0}\cap N$ is cyclic. 
Now let  $N_{r}$ be any Sylow $r$-subgroup of $N$, where $r\in  \pi 
\setminus \{ q \}$. Then for some Sylow $r$-subgroup $G_{r}$ of $G$ we have 
$N_{r}=G_{r}\cap N$   and  $$N\cap G_{r}=(A\cap N)(N\cap 
G_{r})=A_{1}N_{r}=N_{r}A_{1}$$ by Lemma 2.9.
Therefore the hypothesis holds for $(N, A_{1})$.

Finally, in view  of (1) and (4),  $G$ is a $C_{\pi}$-group by Lemma 2.1 (i), which
 contradicts the choice of $G$.

{\bf Corollary  3.2}. {\sl  Let $X$ be a normal  $C_{\pi}$-subgroup of $G$ and 
 $A$ a subgroup of $G$ such     
  that $|G:A|$ is a $\pi $-number. Suppose that 
  $A$ has a  Hall  $\pi$-subgroup $A_{0}$ such that either $A_{0}$ is nilpotent or
 every Sylow subgroup of $A_{0}$ is cyclic.   Suppose that $G$ has a 
subgroup $T$ such that $|G:T|$ is a $\pi'$-number, 
$G=AT$ and 
 $A$  
$X$-permutes  with every Sylow $p$-subgroup of $T$ for all primes
 $p\in  \pi (G)\setminus \pi$,  $p\ne q\in \Bbb{P}$.
 Then $G$ is a $C_{\pi'}$-group.   }
                                            
{\bf Proof.} Let $p\in  \pi (G)\setminus \pi$,   $p\ne q$  and $P$ be 
a Sylow $p$-subgroup of $G$. Since $G=AT$ and  $|G:T|$ is a $\pi'$-number,
 every Sylow $p$-subgroup $T_{p}$ of $T$ 
is a Sylow $p$-subgroup of $G$. Hence for some $x\in G$ we have   
$T_{p}=P^{x}$, so $A$ permutes with $P$ by Lemma 2.10. Therefore Corollary 3.2 follows
 from Proposition 3.1.

{\bf Corollary  3.3} ( see Theorem 1.2 in \cite{GSSI}). {\sl Let $X$ be a normal 
  nilpotent subgroup of $G$ and
  $A$ a Hall $\pi$-subgroup of $G$. 
 Suppose $G=AT$ and 
 $A$  
$X$-permutes  with every subgroup of $T$
 Then $G$ is a $C_{\pi'}$-group.   }

{\bf Corollary 3.4} (See Theorem B in \cite{Cinica}). {\sl Let $A$ be a
 Hall subgroup of a group $G$
and $T$ a minimal supplement of $A$ in $G$. Suppose that  $A$  permutes 
with all Sylow subgroups of $T$ and with all maximal subgroups of
any Sylow subgroup of $T$. Then $G$ is a $C_{\pi'}$-group.  }

{\bf Corollary 3.5} (see Theorem A* in \cite{Proc}). {\sl 
  Let  $H$ be a Hall $\pi$-subgroup of $G$.
Let $G=HT$ for some subgroup $T$ of $G$, and $q$ a prime. If $H$
permutes with every Sylow $p$-subgroup of $T$ for all primes $p\ne
q$, then $T$ contains a complement of $H$ in $G$ and any two
complements of $H$  in $G$ are conjugate.}

{\bf Proof.}  It is clear that $T\cap H$ is a Hall  $\pi$-subgroup of $T$.
Moreover, if $P$ is a Sylow subgroup of $T$ and $HP=PH$, then $HP\cap 
T=P(H\cap T)=(H\cap T)P$. Hence by Proposition 3.1, $T$ is a 
$C_{\pi}$-group. Now   Corollary 3.5 follows from  Proposition 3.1.

{\bf Corollary  3.6} (see \cite{Fog}). {\sl Let   $A$ be a Hall $\pi$-subgroup of $G$,   
  $G=AT$  and  $A$ permutes  with every subgroup of $T$. Then $G$ is an $E_{\pi'}$-group.   }

\section{Proof of Theorem 1.1}

Theorem 1.2 is a special case (when $X=1$) of the following theorem.

{\bf Theorem 4.1.} {\sl Let $X$ be a normal $C_{\pi}$-subgroup
of $G$. Suppose that $G$ has a   subgroup series  $$1=H _0 < H_1 \leq
\ldots \leq H_{t-1}\leq H_t =G,$$  where $|H_{i+1}:H_i|$ is either a $\pi$-number
 or a $\pi'$-number for all $i=1, \ldots , t$.  Suppose that $G$ has a 
 subgroup $T$  such that $H_{1}T=G$  and  $|G:T|$ is  a $\pi'$-number.   
 If $H_{i}$ is $X$-permutable with each  Sylow subgroup of 
$T$,  for    all $i= 1,     \ldots , t$, then $G$ is 
 a $C_{\pi}$-group.}

{\bf Proof.}  Suppose that this theorem  is false and let $G$ be a
counterexample of minimal order.  By Lemma 2.12,  $G$ has a 
Hall $\pi$-subgroup $S$. Hence some Hall $\pi$-subgroup $S_{1}$ of $G$ is not
conjugated with $S$.  Without loss of generality, we may
assume that $H_{1}\ne 1$.  Since $|G:T|$ is        
$\pi'$-number,  every Sylow $p$-subgroup      $P$   of $T$, where $p\in \pi$, is also a
 Sylow $p$-subgroup     of $G$.   
 We proceed the proof via the following
steps.

(1) {\sl $G/N$ is  a $C_{\pi}$-group   for  every non-trivial
quotient   $G/N$ of $G$}.

We consider the subgroup series $$1=H _0 N/N  \leq H_1 N/N \leq \ldots \leq
H_{t-1} N/N \leq H_t N/N=G /N \eqno(2)$$    in $G/N$. Then  $(H_1 N/N)(T 
N/N)=G/N $ and $|G/N:TN/N|=|G:TN|$ is  a $\pi'$-number.   Moreover,  by 
 Lemma 2.6, $H_{i} N/N$ is $XN/N$-permutable with  every Sylow subgroup of $TN/N$
 for all $i= 1,  \ldots , t$.
 On the other
hand, since
$$|H_{i+1}N/N:H_{i}N/N|=|H_{i+1}N:H_{i}N|=|H_{i+1}:H_{i}|:|N\cap
H_{i+1}:N\cap H_{i}|,$$ every index of the  series $ (2)$ is either
$\pi$-number or $\pi'$-number. Moreover, $XN/N \simeq X/(X\cap N)$ is  a $C_{\pi}$-group
by Lemma 2.2 (i). All these show that the hypothesis holds for n $G/N$.
Hence in the case, where $N\ne 1$, $G/N$ is  a $C_{\pi}$-group by the choice of $G$.

(2) $O_{\pi'}(G)=1= O_{\pi}(G).$

Suppose that $D=O_{\pi'}(G)\ne 1$. Then by (1),  there is an  element
$x\in G$ such that ${S_{1}}^{x}D=SD$. But by
 the  Schur-Zassenhaus theorem,
${S_{1}}^{x}$ and $S$ are conjugate in $SD$, which implies that ${S_{1}}$ and $S$
 are conjugate in $G$.
  This contradiction shows
that $O_{\pi'}(G)=1$. Analogously, one can prove that $O_{\pi}(G)=1$.

(3) $X=1$ (This follows from (1), Lemma 2.1 (i)  and the choice of $G$).

(4) $T\neq G$.

Suppose that $T= G$. Then by hypothesis and (3), $H_{1}$
permutes with all Sylow subgroups of $G$. It follows from Lemma 2.7
that $H_{1}$ is subnormal in $G$. Since $H_{1}$ is either a
$\pi$-group or a $\pi'$-group, $H_{1}\leq O_{\pi}(G)$ or $H_{1}\leq
O_{\pi'}(G)$ by Lemma 2.8. It follows from (2) that $H_{1}=1$, which
contradicts to our   assumption about  $H_{1}$. Hence (4) holds.

(5) {\sl The
hypothesis holds for $T$.
}

Consider the  subgroup series  $$1=H _0 \cap T \leq H_1 \cap T \leq
\ldots \leq H_{t-1} \cap T\leq H_t \cap T =T   \eqno(3)$$ of the group 
$T$.  Since
$$H_{i+1}=H_i T \cap H_{i+1}=H_i (H_{i+1}\cap T )$$ we have 
$$|H_{i+1}:H_i|= |H_{i+1}\cap T:H_i \cap T|,$$ for all $ i=1,
 \ldots , t-1$ and $|H_{1}\cap T:H_{0} \cap T|= |H_{1}\cap
T|$ divides    $ |H_1:1|$, we see that every index of the series $(3)$ is
either $\pi$-number or $\pi'$-number.  Now let $E$ be a Sylow subgroup of
$T$. By (3) and the hypothesis, $H_i E=EH_i$. Hence $$H_i E\cap
T =E(H_i \cap T) = (H_i \cap T)E.$$ This shows that the
hypothesis holds for $T$.

(6) {\sl $T$ is  a $C_{\pi}$-group.}

 Since $T\ne G$ by (4), and the
hypothesis holds for $T$ by (5), the minimal choice of $G$ implies that
(6) holds.

(7) {\sl $T$ is  a $E_{\pi'}$-subgroup.} 

This follows from (3), (5) and Lemma 2.12.

(8) {\sl Let $T_{\pi'}$   be a  Hall $\pi'$-subgroup  of $T$   and $D$ 
  a normal subgroup of $G$. Then $T_{\pi'}\ne 1$,  and if
 either  $H_{1}\leq D$ or $T_{\pi'}\leq D$, then  $D=G$}.

Suppose that $T_{\pi'}=1$. Then $H_{1}$ is  a  Hall $\pi'$-subgroup  of 
$G$. Therefore $G$ is a $C_{\pi}$-group  by Corollary 3.2, a contradiction. 
Hence $T_{\pi'}\ne 1$.

We  show that  the  hypothesis holds for $D$. 
 Consider the subgroup series  $$1=D _0\leq  D_1\leq \ldots \leq
D_{t-1}\leq D_t =D, $$ where $D_{i}=H_{i}\cap D$ for all $i=1, 
\ldots , t$. Let $T_{0}=D\cap T$. First we show that  $D_1 T_{
0}=D$. If 
$H_{1}\leq D$, then  $$D=H_{1}(D\cap T)=H_{1} T_{0}=D_1 T_{0}.$$ Now suppose 
that  $T_{\pi'}\leq D$.   In view of 
  (3), $H_{1}$ permutes with $T_{\pi'}$.  Since $H_{1}T=G$, $T\ne G$ and $|G:T|$ is  a 
$\pi'$-number, $H_{1}$ is a  $\pi'$-group. Hence $H_{1}T_{\pi'}$  is a 
Hall $\pi'$-subgroup of $G$. Therefore 
 $$D=(D\cap
H_{1}T_{\pi'})(D\cap T_{\pi})=T_{\pi'}(D\cap H_{1})(D\cap T_{\pi})=(D\cap
H_{1})(T_{\pi'}(D\cap T_{\pi}))=(D\cap H_{1})(D\cap T)=D_1 T_{0}.$$  It is  
also clear that  $|D:T_{0}|$ is  a $\pi'$-number.  Now let $P$ be a Sylow 
 $p$-subgroup of  $T_{0}$.  Then for some Sylow
 $p$-subgroup $T_{p}$  of $T$ we have 
  $P=T_{p}\cap D$. 
 Hence in view of  Lemma  2.9,   $$D\cap H_{i}T_{p}=(D\cap H_{i})(D\cap
 T_{p})=D_{i}P=PD_{i}.$$
    Thus for any $p\in \pi $,   $D_{i}$
 is permutable with every  Sylow $p$-subgroup of 
$T_{0}$  for    all $i= 1,  \ldots , t$.  
 Finally, since the number $$|D_{i}:D_{i-1}|=|(D\cap 
H_{i})H_{i-1}:H_{i-1}|$$ divides $
|H_i:H_{i-1}|$,  each index $|D_{i}:D_{i-1}|$ is either a
$\pi$-number or a $\pi'$-number.  Therefore   the  hypothesis holds for $D$. 
Suppose that $D \ne G$.  Then  $D$ is a $C_{\pi}$-group
 by the choice of $G$.  
 Since  either  $1\ne H_{1}\leq
D$ or $1 \ne T_{\pi'}\leq D$, $G/D$ is a $C_{\pi}$-group  by (1). It follows from
 Lemma 2.1 (i) that $G$ is a
$C_{\pi}$-group, which contradicts the choice
of $G$. Hence, (8) holds.

 {\sl Final contradiction}. 
 Since $G=H_1T$  and  $H_{1}$ permutes with all Sylow subgroups of $T$ by (3),   
$$H_{1}(T_{\pi'})^{x}=(T_{\pi'})^{x}H_{1}$$ for all $x\in G$ by Lemma 2.10.
Therefore  by Lemma 2.3, either $ {H_{1}}^{G}\ne G$ or $ {(T_{\pi'})}^{G}\ne G$. But
in view of (8) both these cases are impossible. 
 The contradiction completes the proof of the result.

{\bf Corollary 4.2  } (see Theorem 5.1 in \cite{Proc}).  {\sl Let $X$ be a normal
 ${\pi}$-separable  subgroup
of $G$. Suppose that $G$ has a   subgroup series  $$1=H _0 < H_1 \leq
\ldots \leq H_{t-1}\leq H_t =G,$$  where $|H_{i+1}:H_i|$ is either a $\pi$-number
 or a $\pi'$-number for all $i=1, \ldots , t$.  Suppose that $G$ has a 
 subgroup $T$  such that $H_{1}T=G$  and  $|G:T|$ is  a $\pi'$-number.   
 If $H_{i}$ is $X$-permutable with each   subgroup of 
$T$,  for    all $i= 1,     \ldots , t$, then $G$ is 
 a $C_{\pi}$-group.}

Theorem 1.1 is a special case (when $X=1$) of the following theorem.

{\bf Theorem 4.3.} {\sl Let $X$ be a normal an $E_{\pi}$-subgroup 
of $G$ and $X_{\pi}$    a  Hall  $\pi$-subgroup of $X$. Suppose that $G$ has a 
  subgroup series  $$1=H _0 < H_1 \leq
\ldots \leq H_{t-1}\leq H_t =G,$$  where $|H_{i+1}:H_i|$ is either a $\pi$-number
 or a $\pi'$-number for all $i=1, \ldots , t$.  Suppose that $G$ has a 
 subgroup $T$  such that $H_{1}T=G$  and  $|G:T|$ is  a $\pi'$-number.   }

(i)    {\sl Suppose that  the Sylow subgroups of $X_{\pi}$  are cyclic. 
 If  $H_{i}$ is $X$-permutable with each  cyclic subgroup  $H$ of 
$T$ of prime power order, for    all $i= 1,  \ldots , t$, then $G$ is 
 a $D_{\pi}$-group.}

(ii)    {\sl Suppose that  $X$ is a $D_{\pi}$-subgroup. 
 If  $H_{i}$ is $X$-permutable with each  cyclic subgroup  $H$ of 
$T$ of prime power order, for    all $i= 1,  \ldots , t$, then $G$ is 
 a $D_{\pi}$-group.}

{\bf Proof.}   (i)    Suppose that this assertion  is false and let $G$ be a
counterexample of minimal order.  
 In view of   Theorem 4.1, $G$ is a $C_{\pi}$-group.  Hence there  is a  $\pi$-subgroup $U$
 of $G$ such that for any  Hall $\pi$-subgroup $E$  of $G$  we have 
$U\nleq E$.

(1) {\sl $G/N$  is a $D_{\pi}$-group for any non-identity normal subgroup $N$ of $G$.}

Consider the subgroup series $$1=H _0 N/N  \leq H_1 N/N \leq \ldots \leq
H_{t-1} N/N \leq H_t N/N=G /N.$$ It is clear that 
$(H_{1}N/N)(TN/N)=G/N$,  $|G/N:TN/N|$ is  a $\pi'$-number 
and 
 $|H_{i+1}N/N:H_i N/N|$ is either a $\pi$-number
 or a $\pi'$-number for all $i=1, \ldots , t$ (see  (1) in the proof of Theorem 4.1).
 Now let $H/N$ be a cyclic subgroup 
of $TN/N$ of prime power order $|H/N|$. Then $H=N(H\cap T)$. Let $W$ be a group of
 minimal order with the 
properties that $W\leq  H\cap T$ and $NW=H$. If $N\cap W\nleq \Phi (W)$, then 
for some maximal subgroup $S$ of $W$ we have $(N\cap W)S=W$. Hence 
$H=NW=N(N\cap W)S=NS$, a contradiction. Hence $N\cap W\leq \Phi (W)$. 
Since $W/W\cap N\simeq    H/N$ is a cyclic group of prime power order, it follows that
 $W$ is cyclic  group of prime power order. Hence $H_{i}$ is 
$X$-permutable with $W$.    Thus  $H_{i}N/N$ is  $XN/N$-permutable with $WN/N=H/N$ by Lemma 2.6.
 Therefore  the hypothesis holds for  $G/N$.
But since $N\ne 1$, $|G/N| < |G|$ and  so $G/N$  is a $D_{\pi}$-group 
 by the choice of $G$.

(2)     $X=1$.

Suppose that $X\ne 1$. Then $G/X$ is  a $D_{\pi}$-group by (1). Hence $G$ 
is  a $D_{\pi}$-group by Lemma 2.2 (iii), a contradiction. Thus we have (2).

(3) {\sl $H_{1}$  permutes with every subgroup of $U$.   }

Let $Z$ be any   cyclic subgroup 
of $U$ of prime power order  $p^{n}$. Then $p\in \pi$ and  $Z\leq G_{p}$ for some Sylow 
$p$-subgroup $G_{p}$ of $G$. Since   $|G:T|$ is  a $\pi'$-number, there is an element 
$x=ht$   such that   $(G_{p})^{x}\leq T$.   
Hence  
$H_{1}Z=ZH_{1}$ by Lemma 2.10, which in view of Lemma 2.13 implies that    
$H_{1}U=UH_{1}$.

(4) $T\neq G$ (see the proof of (4) in the proof of Theorem 4.1).

(5) {\sl $T$  is a  $D_{\pi}$-group. }

The hypothesis holds for $T$ (see (5) in the proof of Theorem 4.1), so in view of (4)
 the minimal 
choice of $G$ implies that we have (5).

(6)  {\sl $V=H_{1}U$ is a  $C_{\pi}$-group.  }

Indeed,  $V$ is a group by (3), and since  $T\neq G$, $H_{1}$ is a 
Hall ${\pi'}$-subgroup of $V$. Therefore we have (6) by (2) and Corollary  
3.2.

(7)  {\sl $ |V:T\cap V|$ is   a $\pi'$-number and  $T\cap V$ is a  $C_{\pi}$-group.  }

First note that $|G:T|=  |V:T\cap V|$ is   a $\pi'$-number.  But   $H_{1}$ is a 
Hall ${\pi'}$-subgroup of $V$. Hence $V=H_{1}(T\cap V)$, which implies 
that $$|V:H_{1}|=|T\cap V:H_{1}\cap T\cap V)|$$ is a $\pi$-number. Hence  
$A= H_{1}\cap T\cap V$  is  a Hall $\pi'$-subgroup of $T\cap V$. Finally, 
if $W$ is any $\pi$-subgroup of $T\cap V=$, then $H_{1}W=WH_{1}$ by (3). 
Therefore  $$AW=(H_{1}\cap T\cap V)W=(H_{1}W\cap T\cap V)=WA.$$   
 Hence 
$T\cap V$ is a  $C_{\pi}$-group by Proposition 3.1.

{\sl Final contradiction for (i).}  In view of (7),  
 $|V:T\cap V|$ is   a $\pi'$-number  and $T\cap V$ is a  $C_{\pi}$-group.  
Hence  in view of (6), 
 there is an element  $x\in V$ such that  $U^{x} \leq    T\cap V\leq T$.  
But by (5), $T$ is  a  $D_{\pi}$-group. Hence for some Hall $\pi$-subgroup 
$T_{\pi}$ of $T$ we have    $U\leq T_{\pi}$, which is a contradiction since  
 $T_{\pi}$ is clearly a  Hall $\pi$-subgroup 
 of $G$. 

(ii)   See  the proof of (i) and use   Lemma 2.2 (iii).

\section{Groups with soluble Hall $\pi$-subgroups}

A group $G$ is said to be:

(a)  an \emph{$E^{S}_{\pi}$-group} (\emph{$E^{N}_{\pi}$-group})
   if $G$ has a   soluble  (a nilpotent, respectively) Hall $\pi$-subgroup;

(b)  a \emph{$C^{S}_{\pi}$-group} (\emph{$C^{N}_{\pi}$-group})
   if $G$ has a   soluble  (nilpotent, respectively) Hall $\pi$-subgroup
 and $G$ is a  $C_{\pi}$-group;

(c)  a \emph{$D^{S}_{\pi}$-group}  
  if $G$ is a  $C^{S}_{\pi}$-group and any
$\pi$-subgroup of $G$ is contained in some Hall $\pi$-subgroup of $G$.

 Our next results are  new criteria for a group to be an 
$E^{S}_{\pi}$-group.

{\bf Lemma 5.1.}
{\sl Let $N$ be a normal $E_{\pi}$-subgroup of $G$  and $N_{\pi}$ a Hall 
$\pi$-subgroup of $N$. Suppose that  $N_{\pi}$ is nilpotent. 
   If   $G/N$ is a 
  $D^{S}_{\pi}$-group, then $G$ is a  $D^{S}_{\pi}$-group }  
 (See \cite{HallIII} or \cite[Chapter IV, Theorem 18.15]{26}).

{\bf Theorem 5.2.} {\sl Let $X$ be a normal $E^{N}_{\pi}$-subroup
 of $G$. Suppose that $G$ has a   subgroup series  $$1=H _0 < H_1 \leq
\ldots \leq H_{t-1}\leq H_t =G,$$  where $|H_{i+1}:H_i|$ is divisible by at most 
 one prime in  $\pi$, 
for all $i=1, \ldots , t$.  Suppose that $G$ has a 
 subgroup $T$  such that $H_{1}T=G$  and  $|G:T|$ is  a $\pi'$-number.     
 If $H_{i}$ is $X$-permutable with each  Sylow subgroup of 
$T$,  for    all $i= 1,  \ldots , t$, then $G$ is 
 an $E^{S}_{\pi}$-group.}

{\bf Proof.}  Suppose that this theorem    is false and let $G$ be a
counterexample of minimal order.    

(1) {\sl Every    non-trivial
quotient   $G/N$ of $G$ is an $E^{S}_{\pi}$-group}
  (See (1) in the proof of Theorem 4.1).

(2) $X=1$.

Suppose that $X\ne 1$.   Then $G/X$ of $G$ is an $E^{ S}_{\pi}$-group  
 by (1). Let $E/X$ be a soluble $\pi$-Hall subgroup of $G/X$. 
Then, by Lemma 5.1, $E$ is  a $D^{S}_{\pi}$-group and if $U$ is a Hall  
subgroup of $E$, then $U$ is a Hall subgroup of $G$. Hence $G$ is
 an $E^{S}_{\pi}$-group,  a contradiction.  Thus we have (2).

(3) $T\neq G$.

Suppose that $T= G$. Then by hypothesis and (2), $H_{1}$
permutes with all Sylow subgroups of $G$.  Hence $H_{1}$  is subnormal in $G$ by Lemma 2.7,
 so
 $$1 < H_{1}\leq 
O_{{\pi'}\cup \{p\}}(G)$$ for some $p\in \pi$ by Lemma 2.8. 
Therefore $O_{{\pi'}\cup \{p\}}(G)$ is a non-identity normal
 $D^{ N}_{\pi}$-group-subgroup of  $G$, which as in the proof of (2), 
conducts us to the contradiction. Thus (3) holds.

{\sl Final contradiction.}
In view of (2), the hypothesis is true for $T$ (see (5) in the proof of Theorem 4.1). 
Hence   $T$ is an   $E^{S}_{\pi}$-group by (3) and the choice of 
$G$.  But since $|G:T|$ is  a $\pi'$-number, any Hall $\pi$-subgroup of $T$ is a Hall 
$\pi$-subgroup of $G$ as well.  Hence   $G$ is an   
$E^{S}_{\pi}$-group,  a contradiction. 

The theorem is proved.

{\bf Theorem 5.3.} {\sl  Let $X$ be a normal $E^{N}_{\pi}$-subgroup
 of $G$. Suppose that $G$ has a   subgroup series  $$1=H _0 < H_1 \leq
\ldots \leq H_{t-1}\leq H_t =G,$$  where $|H_{i+1}:H_i|$ is divisible by at most 
 one prime in  $\pi$, 
for all $i=1, \ldots , t$.  Suppose that $G$ has a 
 subgroup $T$  such that $H_{1}T=G$  and  $|G:T|$ is  a $\pi'$-number.     
  If  $H_{i}$ is $X$-permutable with each  cyclic subgroup  $H$ of 
$T$ of prime power order, for    all $i= 1,  \ldots , t$, then $G$ is 
 a $D^{S}_{\pi}$-group.}

{\bf Proof.}   Suppose that this theorem    is false and let $G$ be a
counterexample of minimal order.   Then $|\pi \cap \pi (G)| > 1$. Without
 loss we may suppose that  $H_{1}\ne 1$  and    $H_{t-1}\ne G$.

(1) {\sl Every    on-trivial
quotient   $G/N$ of $G$ is a $D^{ S}_{\pi}$-group}
  (See (1) in the proof of Theorem 4.3).

(2) {\sl If $N$ is a normal $E^{N}_{\pi}$-subgroup of $G$,
 then $N=1$. In particular,  $X=1$. }

Suppose that $N\ne 1$  Then $G/N$ of $G$ is a $D^{ S}_{\pi}$-group  
by (1) and the choice of $G$. 
 Hence $G$ is  a $D^{ 
S}_{\pi}$-group  by Lemma 5.1, which contradicts the choice of $G$. 
Hence we have (2).

(3) $T\neq G$ (See (3) in the proof of Theorem 5.2).

(4) {\sl The hypothesis holds for  $H_{t-1}$  and  $T$.}

Since $$H_{t-1}= H_{1}(H_{t-1}\cap T)$$  and   $$|H_{t-1}:H_{t-1}\cap 
T|=|G:T|,$$ the hypothesis holds for  $H_{t-1}$ by (2). The second assertion of (4) may
 be proved as (5) in the proof of Theorem 4.1.

(5) {\sl  $G$ is a $C^{ S}_{\pi}$-group}.

Since  $|G:H_{t-1}|$ is divisible by at most 
 one prime in  $\pi$ and  $|\pi \cap \pi (G)| > 1$, there is a prime
 $p\in \pi \cap \pi (G)$ such that  for a Sylow $p$-subgroup $P$ of $G$ we have
 $P\leq H_{t-1}$, so in view of (2) and Lemma 2.10  we have $1 < P \leq (H_{t-1})_{G}$.  
Since $H_{t-1}\ne G$ and the hypothesis holds for $H_{t-1}$,  $H_{t-1}$ is a
 $D^{ S}_{\pi}$-group by the choice
 of $G$.  Hence $(H_{t-1})_{G}$ is  a $C^{ S}_{\pi}$-group, so $G$ is a
 $C^{ S}_{\pi}$-group by Lemma 2.1 (i) and Theorem 5.2. 

 (6) {\sl   For any $i$, $H_{i}$ permutes with every $\pi$-subgroup of 
$G$}  (see (3) in the proof of Theorem 4.3).

(7) {\sl     $V=H_{1}U$  is a $D^{\cal S}_{\pi}$-group for any
 $\pi$-subgroup $U$  of $G$.}

 By (6), $V$ is a subgroup of $G$. Moreover, $V$ is a $C_{\pi}$-group by Proposition  3.1.
 We show that the hypothesis holds for $V$. 
Let $V_{\pi}$ be a  Hall $\pi$-subgroup of $V$. Then $V=H_{1}V_{\pi}$
  and  $$V_{i}=V\cap H_{i}=H_{1}(V_{\pi}\cap H_{i}),$$  for all $i=1, \ldots , t$. 
Let $W$ be any subgroup of  $V_{\pi}$. Then
 $$V_{i}W=(H_{1}(V_{\pi}\cap H_{i}))W=H_{1}(V_{\pi}\cap H_{i}W)=
H_{1}W(V_{\pi}\cap H_{i})W(H_{1}(V_{\pi}\cap H_{i}))=WV_{i}.$$ It is clear also that
 $|V_{i+1}:V_i|$ is divisible by at most 
 one prime in  $\pi$, 
for all $i=0, \ldots , t$. Hence the hypothesis holds for $V$.   Suppose 
that $G=V$. In this case, in view of (5), we may suppose that $V_{\pi}$ is 
a soluble Hall $\pi$-subgroup of $G$. Let $L$ be a minimal normal subgroup 
of $V_{\pi}$. Then $$L^{G}=L^{V_{\pi}H_{1}}=L^{H_{1}}\leq 
 LH_{1}\cap L^{H_{1}}=L(L^{G}\cap H_{1})$$
 and $L$ is
 a $q$-group for some $q\in \pi$. From (2) it follows that for some prime 
$p\in \pi$ with $q\ne p$ 
 we have $p\in 
\pi (L^{G}\cap H_{1})$ and a Sylow $p$-subgroup $P$ of $L^{G}\cap H_{1}$ is also a
 Sylow 
subgroup of $L^{G}$. For some Sylow $p$-subgroup  $G_{p}$ of $G$ we have  
$P=L^{G}\cap  G_{p}$, so $$L^{G}\cap  H_{1}G_{p}=(L^{G}\cap  H_{1})(L^{G}\cap G_{p})=(L^{G}\cap
  H_{1})P=P(L^{G}\cap  H_{1}).$$ 
 Hence $P\leq (L^{G}\cap H_{1})_{L^{G}}.$ Therefore $G$ has a
 non-identity
 subnormal subgroup  $R=(L^{G}\cap H_{1})_{L^{G}}$ of order divisible by 
at most one prime in $\pi$, which contradicts (2).  Hence $V\ne G$, so  
$V$  is a $D^{ S}_{\pi}$-group  by the choice of $G$.

{\sl Final contradiction.}  Let $U$ be any $\pi$-subgroup of $G$ and $V=H_{1}U$. 
Then $V$ is a $D^{ S}_{\pi}$-subgroup of $G$ by (7), and
 $|G:T|=|V:V\cap T|$ is a $\pi'$-number.  
 By (3) and (4), $T$ is   a 
$D_{\pi}$-group.   Hence for some Hall 
$\pi$-subgroup $V_{\pi}$  of $V$ we have $U\leq V_{\pi}$,  and $ 
V_{\pi}\leq  T_{\pi}$ . But since $|G:T|$ is a $\pi'$-number, 
 $T_{\pi}$ is a Hall $\pi$-subgroup   of $G$. Therefore
$G$ is a $D^{ S}_{\pi}$-group.

 {\bf Corollary  5.4.} {\sl Suppose that $G$ has a   subgroup series  $$1=H _0 < H_1 \leq
\ldots \leq H_{t-1}\leq H_t =G,$$  where $|H_{i+1}:H_i|$ is divisible by at most 
 one prime in  $\pi$,   for all $i=1, \ldots t.$   Suppose also that $G$ has a 
 subgroup $T$  such that $H_{1}T=G$  and  $|G:T|$ is  a $\pi'$-number.   If  $H_{i}$
 is quasinormal in $T$,  for    all $i= 1,  \ldots , t$,
 then $G$ is  a $D^{S}_{\pi}$-group.}

{\bf Corollary 5.6 } (see S. A. $\breve{C}$unihin \cite{Pi-otdel} or \cite[Chapter IV,
 Theorem   18.13]{26}). {\sl  If $G$  has a  chief 
 subgroup series  $$1=H _0 < H_1 \leq
\ldots \leq H_{t-1}\leq H_t =G,$$  where $|H_{i+1}:H_i|$ is divisible by at most 
 one prime in  $\pi$, 
for all $i=1, \ldots , t$, then $G$ is 
 a $D^{ S}_{\pi}$-group.}

\end{document}